**Hankel determinants of middle binomial coefficients and conjectures for some polynomial extensions and modifications.**

Johann Cigler
johann.cigler@univie.ac.at

**Abstract**
The middle binomial coefficients can be interpreted as numbers of Motzkin paths which have no horizontal steps at positive heights. Assigning suitable weights gives some nice polynomial extensions. We determine the Hankel determinants and their generating functions for the middle binomial coefficients and derive many conjectures for their polynomial extensions. Finally, we explore experimentally some modifications of the middle binomial coefficients whose Hankel determinants show an interesting modular pattern and obtain some q-analogs.

## 1. Introduction

Following [2], we shall call the numbers

$$b(n)=\binom{n}{\left\lfloor \frac{n}{2}\right\rfloor} \tag{1}$$

**middle binomial coefficients** in order to distinguish them from the central binomial coefficients $b(2n)=\binom{2n}{n}$. In OEIS [10], A001405, some combinatorial interpretations of these numbers are mentioned. For our purposes we interpret $b(n)$ as the number of Motzkin paths of length $n$ which have no horizontal steps at positive heights. Let

$$D_k(n)=\det\left(b(k+i+j)\right)_{i,j=0}^{n-1} \tag{2}$$

denote their Hankel determinants. The first two Hankel determinants are well known (cf. [2]). The original aim of the present paper was the determination of all Hankel determinants $D_k(n)$ and their generating functions. Thereby we were led to explore in an experimental way some polynomial extensions by assigning various weights depending on a parameter $t$. Finally, we state some conjectures about the Hankel determinants $D_k^{(r)}(n)=\det\left(b^{(r)}(k+i+j)\right)_{i,j=0}^{n-1}$ of the sequences $b^{(r)}(n)=\binom{n}{\left\lfloor \frac{n-r}{2}\right\rfloor}$ for $r\in\mathbb{N}$, which show interesting modular patterns and compute the Hankel determinants of the $q-$analogs $b(n,q)=\begin{bmatrix} n \\ \left\lfloor \frac{n}{2}\right\rfloor\end{bmatrix}_q \frac{2^n}{(-q;q)_{\left\lfloor \frac{n}{2}\right\rfloor}(-q;q)_{\left\lfloor \frac{n+1}{2}\right\rfloor}}.$

I want to thank Max Alekseyev, Sam Hopkins and Christian Krattenthaler for valuable hints and proofs and N.J.A. Sloane for founding The On-Line Encyclopedia of Integer Sequences which provided a wealth of useful information.

## 2. Some nice polynomial extensions of the middle binomial coefficients

**2.1.** Let $\mathbf{B}_n$ denote the set of all Motzkin paths of length $n$ with no horizontal steps at positive heights. $\mathbf{B}_n$ consists of lattice paths in $\mathbb{N}\times\mathbb{Z}$ from $(0,0)$ to $(0,n)$ which remain at or above the $x-$axis with up-steps $U_j:(i-1,j-1)\to(i,j)$, down-steps $D_j:(i-1,j+1)\to(i,j)$ and horizontal steps $H_i:(i-1,0)\to(i,0)$.

Let us first verify that $\mathbf{B}_n$ has $b(n)=\binom{n}{\left\lfloor \frac{n}{2}\right\rfloor}$ elements. To this end we show more generally that the number $b(n,k)$ of all such paths from $(0,0)$ to $(n,k)$ with no horizontal steps at positive heights is $b(n,k)=\binom{n}{\left\lfloor \frac{n-k}{2}\right\rfloor}$.

The numbers $b(n,k)$ are uniquely determined by the following properties:
$b(0,k)=[k=0]$,
$b(n,0)=b(n-1,0)+b(n-1,1)$,
$b(n,k)=b(n-1,k-1)+b(n-1,k+1)$ for $k>0$.

Since the numbers $c(n,k)=\binom{n}{\left\lfloor \frac{n-k}{2}\right\rfloor}$ satisfy the same properties they coincide with $b(n,k)$.

If we assign to all down-steps of $\mathbf{B}_n$ the weight $t$ and to all other steps the weight 1 then the weight becomes

$$a_n(t)=\sum_{j=0}^{\left\lfloor \frac{n}{2}\right\rfloor}\left(\binom{n}{j}-\binom{n}{j-1}\right)t^j. \tag{3}$$

To prove (3) it suffices to show that the set $\mathbf{B}_{n,k}$ of all elements of $\mathbf{B}_n$ with $k$ down-steps has $\binom{n}{k}-\binom{n}{k-1}$ elements. As indicated in OEIS [10], A008315, $\mathbf{B}_{n,k}$ can bijectively be mapped onto the set $\mathbf{A}_{n,k}$ of all paths from $(0,0)$ to $(n,n-2k)$ consisting of up-steps and down-steps which never go below the $x-$axis by replacing the $n-2k$ horizontal steps with up-steps. The reflection principle then shows that $\left|\mathbf{A}_{n,k}\right|=\binom{n}{k}-\binom{n}{k-1}$. For there are $\binom{n}{k}$ paths consisting only of up-steps and $k$ down-steps without restrictions. Consider those paths which go below the $x-$axis. Let $D_{-1}$ be the first step with this property. If we reflect the remaining path on the line $y=-1$ we get a a bijection to all $\binom{n}{k-1}$ paths from $(0,0)$ to $(n,2k-n-2)$.

**2.2.** By assigning the weight $t$ only to the down-steps $D_{2j}$, $j \geq 0,$ and to all other steps the weight 1 then the weight of $\mathbf{B}_n$ becomes

$$b_n(t) = \sum_{j=0}^{n} \binom{\left\lfloor \frac{n}{2} \right\rfloor}{j} \binom{\left\lfloor \frac{n+1}{2} \right\rfloor}{j} t^j. \tag{4}$$

It suffices to verify that the polynomials

$$\begin{aligned} c(n,2k,t) &= \sum_{j=0}^{n} \binom{\left\lfloor \frac{n}{2} \right\rfloor}{j+k} \binom{\left\lfloor \frac{n+1}{2} \right\rfloor}{j} t^j, \\ c(n,2k-1,t) &= \sum_{j=0}^{n} \binom{\left\lfloor \frac{n}{2} \right\rfloor}{j} \binom{\left\lfloor \frac{n+1}{2} \right\rfloor}{j+k} t^j \end{aligned} \tag{5}$$

satisfy
$c(0,k,t) = [k=0],$
$c(n,0,t) = c(n-1,0,t) + tc(n-1,1,t),$
$c(n,k,t) = c(n-1,k-1,t) + t_k c(n-1,k+1,t),$
where $t_k$ satisfies $t_{2k} = t$ *and* $t_{2k+1} = 1$ which is easily verified.

**2.3.** Consider the set of all Motzkin paths of length $n$ and assign to all up-steps weight 1, to all down-steps weight $t$, to all horizontal steps on height $0$ weight $1,$ to all horizontal steps on even heights weight $1-t$ and to all horizontal steps on odd heights weight $t-1.$ Then the total weight becomes

$$c_n(t) = \sum_{j=0}^{n} \binom{\left\lfloor \frac{n-1}{2} \right\rfloor}{\left\lfloor \frac{j-1}{2} \right\rfloor} \binom{\left\lfloor \frac{n}{2} \right\rfloor}{\left\lfloor \frac{j}{2} \right\rfloor} t^{j-1} \tag{6}$$

with $c_0(t) = 1.$

By Flajolet's continued fraction theorem (cf. [5], V.4) the generating function of these paths is

$$F(z) := \sum_{n \geq 0} c_n(t) z^n$$

$$= \cfrac{1}{1 - z - \cfrac{tz^2}{1-(t-1)z - \cfrac{tz^2}{1-(1-t)z - \cfrac{tz^2}{1-(t-1)z - \cfrac{tz^2}{1-(1-t)z - \ddots}}}}}.$$

With

$$G(z):=\cfrac{1}{1-(t-1)z-\cfrac{tz^2}{1-(1-t)z-\cfrac{tz^2}{1-(t-1)z-\cfrac{tz^2}{1-(1-t)z-\cfrac{tz^2}{1-(t-1)z-\ddots}}}}}$$

we get

$$F(z)=\frac{1}{1-z-tz^2G(z)}. \tag{7}$$

$G(z)$ satisfies the identity

$$G(z)=\cfrac{1}{1-(t-1)z-\cfrac{tz^2}{1-(1-t)z-tz^2G(z)}}.$$

Solving this quadratic equation gives

$$G(z)=\frac{-t^2z^2+2tz^2-z^2+1-\sqrt{t^4z^4-2t^2z^4-2t^2z^2+z^4-2z^2+1}}{2tz^2\left(-tz+z+1\right)}.$$

Therefore we get

$$F(z)=\frac{1-2tz-z^2+t^2z^2-\sqrt{1-2z^2\left(1+t^2\right)+z^4\left(1-t^2\right)^2}}{2tz\left(1+(1-t)z\right)}. \tag{8}$$

It remains to show that $\left[z^n\right]F(z)=c_n(t).$ A direct proof is rather tedious. But there is an elegant proof due to Max Alekseyev [3]. Observe that

$$\left[z^{2n+1}\right]z\left(\left(1+z^2\right)\left(1+t^2z^2\right)\right)^n=\left[z^{2n+1}\right]\sum_{\ell+j=n}\binom{n}{\ell}z^{2\ell+1}\binom{n}{j}t^{2j}z^{2j}=\sum_{j=0}^{n}\binom{n}{j}^2t^{2j},$$

$$\left[z^{2n+1}\right]z^3t\left(\left(1+z^2\right)\left(1+t^2z^2\right)\right)^n=\left[z^{2n+1}\right]\sum_{\ell+j=n}\binom{n}{\ell}z^{2\ell+1}t\binom{n}{j-1}t^{2j-2}z^{2j}=\sum_{j=0}^{n}\binom{n}{j-1}\binom{n}{j}t^{2j-1}$$

and $\left[z^{2n}\right]z^jt\left(\left(1+z^2\right)\left(1+t^2z^2\right)\right)^n=0$ for odd $j.$

Therefore

$$\left[z^{2n+1}\right]z\left(1+z^2t\right)\left(\left(1+z^2\right)\left(1+t^2z^2\right)\right)^n=c_{2n+1}(t).$$

In the same way we get

$$\left[z^{2n}\right]\left(1+z^2t\right)\left(1+z^2\right)^n\left(1+t^2z^2\right)^{n-1}=c_{2n}(t).$$

Combining these identities, we get

$$c_n(t)=\left[z^n\right]\left(\left(1+z^2\right)\left(1+z^2t^2\right)\right)^{\frac{n}{2}}\left(1+tz^2\right)\left(z\left(\left(1+z^2\right)\left(1+z^2t^2\right)\right)^{-\frac{1}{2}}+\left(1+t^2z^2\right)^{-1}\right). \tag{9}$$

Now we use the Lagrange-Bürmann formula in the following form:

Let $f(w)=\dfrac{w}{\varphi(w)}$ with $\varphi(0)\neq 0$ and let $g(z)$ be the inverse satisfying $f(g(z))=z.$

Then

$$\left[z^n\right]H(g(z))=\left[w^n\right]H(w)\varphi(w)^{n-1}\left(\varphi(w)-w\varphi'(w)\right). \tag{10}$$

Let $\varphi(z)=\sqrt{\left(1+z^2\right)\left(1+t^2z^2\right)}.$ Then we have $\varphi(w)-w\varphi'(w)=\dfrac{\left(1-w^2t\right)\left(1+w^2t\right)}{\varphi(w)}$ and

$$g(z)=\sqrt{\frac{1-\left(1+t^2\right)z^2-\sqrt{1-2\left(1+t^2\right)z^2+\left(1-t^2\right)^2z^4}}{2t^2z^2}}.$$

With these notations we get

$$\left(\left(1+w^2\right)\left(1+w^2t^2\right)\right)^{\frac{n}{2}}\left(1+tw^2\right)\left(w\left(\left(1+w^2\right)\left(1+w^2t^2\right)\right)^{-\frac{1}{2}}+\left(1+t^2w^2\right)^{-1}\right)$$
$$=\left(1+tw^2\right)\left(w\varphi(w)^{n-1}+\varphi(w)^{n-2}\left(1+w^2\right)\right).$$

We want to choose $H(w)$ such that

$$\left(1+tw^2\right)\left(w\varphi(w)^{n-1}+\varphi(w)^{n-2}\left(1+w^2\right)\right)=H(w)\varphi(w)^{n-1}\left(\varphi(w)-w\varphi'(w)\right)$$
$$=H(w)\varphi(w)^{n-2}\left(1-w^2t\right)\left(1+w^2t\right).$$

This gives

$$H(w)=\frac{w\varphi(w)+1+w^2}{1-tw^2}.$$

Then a computation with Mathematica gives $H(g(z))=F(z)$ and therefore (8) and (9) give (6).

**Remark**

It is perhaps interesting to note that $c_n(t)=C_n(t,-1)$ is the restriction to $q=-1$ of the $q-$ Narayana polynomials

$$C_n(t,q)=\sum_{j=0}^{n}\begin{bmatrix}n\\j\end{bmatrix}_q\begin{bmatrix}n-1\\j\end{bmatrix}_q\frac{q^{j^2+j}}{\left[j+1\right]_q}t^j. \tag{11}$$

For $t=1$ this gives (cf. [6], (3.4)) a $q-$ analog of the Catalan numbers

$$C_n(1,q)=\sum_{j=0}^{n}\begin{bmatrix} n \\ j \end{bmatrix}_q \begin{bmatrix} n-1 \\ j \end{bmatrix}_q \frac{q^{j^2+j}}{[j+1]_q}=\frac{1}{[n+1]_q}\begin{bmatrix} 2n \\ n \end{bmatrix}_q. \tag{12}$$

Note that

$$\lim_{q\to -1}\frac{1}{[n+1]_q}\begin{bmatrix} 2n \\ n \end{bmatrix}_q = \binom{n}{\left\lfloor \frac{n}{2} \right\rfloor}. \tag{13}$$

Another relation with Catalan numbers is $c_{2n+1}(-1)=C_n$.

**3. The Hankel determinants of the middle binomial coefficients and their generating functions.**

Let us first recall some well-known results about orthogonal polynomials (cf. e.g. [4] where they are stated in a form which will be needed here).

A sequence of monic polynomials $p_n(x)$ of degree $n$ is called orthogonal if

$$p_n(x)=(x-s_{n-1})p_{n-1}(x)-t_{n-2}p_{n-2}(x) \tag{14}$$

for some values $s_n$ and $t_n$. Let $\Lambda$ denote the linear functional on the polynomials defined by $\Lambda\left(p_n(x)\right)=[n=0]$. Then (14) implies the orthogonality relations $\Lambda\left(p_n(x)p_m(x)\right)=0$ for $n\neq m$ and $\Lambda\left(p_n(x)^2\right)\neq 0$.

Define numbers $c(n,k)$ by

$$x^n=\sum_{k=0}^{n}c(n,k)p_k(x). \tag{15}$$

They satisfy

$$\begin{aligned} &c(0,k)=[k=0],\\ &c(n,0)=s_0c(n-1,0)+t_0c(n-1,1),\\ &c(n,k)=c(n-1,k-1)+s_kc(n-1,k)+t_kc(n-1,k+1). \end{aligned} \tag{16}$$

**Remark**

Recall that $c(n,k)$ can be interpreted as the weight of the Motzkin paths from $(0,0)$ to $(n,k)$, where the horizontal steps at height $k$ have weight $s_k$ and the down-steps $D_k$ have weight $t_k$. All up-steps have weight 1.

The moments $M(n)=\Lambda\left(x^n\right)$ are given by

$$M(n)=\Lambda\left(x^n\right)=c(n,0). \tag{17}$$

The first two Hankel determinants are

$$\det\left(M(i+j)\right)_{i,j=0}^{n-1} = \prod_{i=1}^{n-1}\prod_{j=0}^{i-1} t_j,$$

$$\det\left(M(i+j+1)\right)_{i,j=0}^{n-1} = (-1)^n p_n(0)\prod_{i=1}^{n-1}\prod_{j=0}^{i-1} t_j. \tag{18}$$

For the middle binomial coefficients, we get $t_j = 1$, $s_0 = 1$ and $s_j = 0$ for $j > 0$.
By (18) we get

$$D_0(n) = \det\left(b(i+j)\right)_{i,j=0}^{n-1} = \prod_{i=1}^{n-1}\prod_{j=0}^{i-1} 1 = 1. \tag{19}$$

The polynomials $p_n(x)$ satisfy $p_0(x) = 1,$ $p_1(x) = x-1$ and $p_n(x) = xp_{n-1}(x) - p_{n-2}(x).$

Let $F_n(x)$ denote the Fibonacci polynomials defined by $F_n(x) = xF_{n-1}(x) - F_{n-2}(x)$ with $F_0(x) = 1$ and $F_{-1}(x) = 0.$ Then

$$p_n(x) = F_n(x) - F_{n-1}(x). \tag{20}$$

Thus $p_n(0) = (-1)^{\binom{n+1}{2}}$ and therefore

$$D_1(n) = (-1)^{\binom{n+1}{2}-n} = (-1)^{\binom{n}{2}}. \tag{21}$$

For the computation of higher Hankel determinants, we use Dodgson's condensation method (cf. [8], Prop. 10) which gives

$$D_k(n)D_{k+2}(n-2) - D_{k+2}(n-1)D_k(n-1) + D_{k+1}(n-1)^2 = 0. \tag{22}$$

**Theorem 1**

$$(-1)^{k\binom{n}{2}} D_k(n) = \prod_{j=1}^{k}\left(\frac{j+n}{j}\right)^{\min(j,k-j)} = \prod_{i=1}^{\left\lfloor \frac{k}{2}\right\rfloor}\prod_{j=1}^{\left\lceil \frac{k}{2}\right\rceil}\frac{i+j+n-1}{i+j-1}. \tag{23}$$

Let us first derive some properties of these numbers.

**Proposition 2**

*Let*

$$r_k(x) = \prod_{j=1}^{k}\left(\frac{j+x}{j}\right)^{\min(j,k-j)},$$

$$s_k(x) = \prod_{i=1}^{\left\lfloor \frac{k}{2}\right\rfloor}\prod_{j=1}^{\left\lceil \frac{k}{2}\right\rceil}\frac{i+j+x-1}{i+j-1}. \tag{24}$$

*Then*

$$s_k(x) = r_k(x) = r_{k-1}(x) \prod_{j=\left\lfloor \frac{k+1}{2} \right\rfloor}^{k-1} \frac{x+j}{j} \tag{25}$$

*and*

$$r_k(x) = \det\left( \binom{x + \left\lfloor \frac{k}{2} \right\rfloor + i + j}{\left\lfloor \frac{k}{2} \right\rfloor + j} \right)_{i,j=0}^{\left\lfloor \frac{k-1}{2} \right\rfloor}. \tag{26}$$

**Proof**

$$s_{2k}(x) = \prod_{i=0}^{k-1} \prod_{j=1}^{k} \frac{i+j+x}{i+j} = \left(\frac{1+x}{1}\right)\left(\frac{2+x}{2}\right)^2 \cdots \left(\frac{k+x}{k}\right)^k \left(\frac{k+1+x}{k+1}\right)^{k-1} \cdots \left(\frac{2k-1+x}{2k-1}\right) = r_{2k}(x),$$

$$s_{2k+1}(x) = \prod_{i=0}^{k-1} \prod_{j=1}^{k+1} \frac{i+j+x}{i+j} = \prod_{i=0}^{k-1} \prod_{j=1}^{k+1} \frac{i+j+x}{i+j} = \left(\frac{1+x}{1}\right)\left(\frac{2+x}{2}\right)^2 \cdots \left(\frac{k+x}{k}\right)^k \left(\frac{k+1+x}{k+1}\right)^{k} \cdots \left(\frac{2k+x}{2k}\right) = r_{2k+1}(x).$$

$$\frac{s_{2k}(x)}{s_{2k-1}(x)} = \frac{\prod_{i=1}^{k} \prod_{j=1}^{k} \frac{i+j+x-1}{i+j-1}}{\prod_{i=0}^{k-1} \prod_{j=1}^{k} \frac{i+j+x-1}{i+j-1}} = \prod_{j=1}^{k} \frac{k+j+x-1}{k+j-1} = \prod_{j=k}^{2k-1} \frac{x+j}{j},$$

$$\frac{s_{2k+1}(x)}{s_{2k}(x)} = \frac{\prod_{i=1}^{k} \prod_{j=1}^{k+1} \frac{i+j+x-1}{i+j-1}}{\prod_{i=0}^{k} \prod_{j=1}^{k} \frac{i+j+x-1}{i+j-1}} = \prod_{i=1}^{k} \frac{k+j+x}{k+j} = \prod_{j=k+1}^{2k} \frac{x+j}{j}.$$

We can also write $s_{2k}(x) = \prod_{j=0}^{k-1} \frac{(x+k+j)!\,j!}{(x+j)!(j+k)!}.$

Now we have

$$\det\left( \binom{x+k+i+j}{k+j} \right)_{i,j=0}^{k-1} = \det\left( \frac{(x+k+i+j)!}{(k+j)!(x+i)!} \right)_{i,j=0}^{k-1}$$

$$= \prod_{j=0}^{k-1} \frac{1}{(k+j)!} \frac{(x+k)!}{x!} \frac{(x+k+1)!}{(x+1)!} \cdots \frac{(x+2k-1)!}{(x+k-1)!} \det\left( \prod_{\ell=1}^{j} (x+k+i+\ell) \right)_{i,j=0}^{k-1}$$

$$=\prod_{j=0}^{k-1}\frac{(x+k+j)!}{(k+j)!(x+j)!}\det\left(\prod_{\ell=1}^{j}(x+k+i+\ell)\right)_{i,j=0}^{k-1}=s_{2k}(x)\prod_{j=0}^{k-1}\frac{1}{j!}\det\left(\prod_{\ell=1}^{j}(x+k+i+\ell)\right)_{i,j=0}^{k-1}$$

$$=r_{2k}(x)\prod_{j=0}^{k-1}\frac{1}{j!}\det\left(p_j(x+i)\right)_{i,j=0}^{k-1}$$

with $p_j(x)=\prod_{\ell=1}^{j}(x+k+\ell).$

By [8], Proposition 1 we get $\det\left(p_j(x+i)\right)_{i,j=0}^{k-1}=\prod_{j=1}^{k-1}\prod_{i=0}^{j}(j-i)=\prod_{j=1}^{k-1}j!.$

Therefore,

$$\det\left(\binom{x+k+i+j}{k+j}\right)_{i,j=0}^{k-1}=r_{2k}(x).$$

The same argument gives

$$\det\left(\binom{x+k+i+j}{k+j}\right)_{i,j=0}^{k}=r_{2k+1}(x).$$

**Proof of Theorem 1**

By (22) it suffices to show that
$r_k(x)$ satisfies

$$r_k(x)r_{k-2}(x)=(-1)^k r_{k-2}(x+1)r_k(x-1)+r_{k-1}(x)^2. \tag{27}$$

We get
$xr_{2k}(x)r_{2k-2}(x)=r_{2k}(x-1)r_{2k-2}(x+1)\left(x+2k-1\right)$
because the denominators are the same and the powers of $(x+j)$ are the same on both sides. This is true for $j=0$ and $j=2k-1.$ For $1\le j\le 2k-2$ the power of $x+j$ on the left -hand side is $\min(j,2k-j)+\min(j,2k-2-j).$ On the right-hand side we get the same value since $\min(j+1,2k-j-1)+\min(j-1,2k-1-j)=1+\min(j,2k-j-2)-1+\min(j,2k-j).$

Further we get

$(2k-1)r_{2k-2}(x)r_{2k}(x)=r_{2k-1}(x)^2\left(x+2k-1\right)$
because
$\frac{s_{2k}(x)}{s_{2k-1}(x)}=\prod_{j=k}^{2k-1}\frac{x+j}{j}$ and $\frac{s_{2k-1}(x)}{s_{2k-2}(x)}=\prod_{j=k}^{2k-2}\frac{x+j}{j}.$

This implies

$$r_{2k}(x)r_{2k-2}(x)-r_{2k-2}(x+1)r_{2k}(x-1)-r_{2k-1}(x)^2=r_{2k}(x)r_{2k-2}(x)\left(1-\frac{x}{x+2k-1}-\frac{2k-1}{x+2k-1}\right)=0.$$

In the odd case analogous arguments give

$xr_{2k+1}(x)r_{2k-1}(x)=r_{2k+1}(x-1)r_{2k-1}(x+1)\left(x+2k\right)$ and $2(x+k)r_{2k+1}(x)r_{2k-1}(x)=r_{2k}(x)^2\left(x+2k\right).$

This implies

$$r_{2k+1}(x)r_{2k-1}(x)+r_{2k-1}(x+1)r_{2k+1}(x-1)-r_{2k}(x)^2=r_{2k+1}(x)r_{2k-1}(x)\left(1+\frac{x}{x+2k}-\frac{2(x+k)}{x+2k}\right)=0.$$

Since the binomial coefficients are integers for $k\in\mathbb{N}$ we see again that the polynomials $r_k(x)$ are integer valued, i.e. $r_k(n)\in\mathbb{N}$ for $n\in\mathbb{N}.$

The first terms of the sequence $\left(r_k(x)\right)_{k\geq 0}$ are $1, 1, 1+x, \frac{1}{2}(1+x)(2+x),$
$\frac{1}{12}(1+x)(2+x)^2(3+x), \frac{1}{144}(1+x)(2+x)^2(3+x)^2(4+x).$

**Remark**

The numbers $\prod_{i=1}^{\left\lfloor\frac{k}{2}\right\rfloor}\prod_{j=1}^{\left\lceil\frac{k}{2}\right\rceil}\frac{i+j+n-1}{i+j-1}$ also count plane partitions contained in a $\left\lfloor\frac{k}{2}\right\rfloor\times\left\lceil\frac{k}{2}\right\rceil\times n-$ box (cf. [7] and [11], 2, 7.21) and occur in many other situations in a natural way, cf. e.g. OEIS [10], A006542, A107915, A140902, A002415, A047819.

Since $r_{2k}(x)$ is a polynomial of degree $k^2$ and $r_{2k+1}(x)$ a polynomial of degree $k^2+k$ we get

$$\sum_{n\geq 0}D_{2k}(n)x^n=\sum_{n\geq 0}r_{2k}(n)x^n=\frac{A_{2k}(x)}{(1-x)^{k^2+1}} \tag{28}$$

and

$$\sum_{n\geq 0}r_{2k+1}(n)x^n=\frac{A_{2k+1}(x)}{(1-x)^{k^2+k+1}} \tag{29}$$

for some polynomials $A_k(x).$

These turn out to be palindromic, unimodal and gamma-nonnegative with degree $\deg A_{2k}(x)=(k-1)^2$ for $k\geq 1$ and $\deg A_{2k+1}(x)=k^2-k.$

Recall that a polynomial $a_n(x)=\sum_{j=0}^{n}a(j)x^j$ is called palindromic if $x^na_n\left(\frac{1}{x}\right)=a_n(x),$ i.e. if its coefficients are the same when read from left to right as from right to left. It is unimodal if its coefficients first weakly increase and then weakly decrease. It is called gamma-nonnegative if it can be written as $\sum_{j=0}^{\left\lfloor\frac{n}{2}\right\rfloor}\gamma_jx^j(1+x)^{n-2j}$ with $\gamma_j\geq 0.$ A gamma-nonnegative polynomial is both unimodal and palindromic.

Let us consider the polynomials $A_k(x)$ for small $k.$ Computations give

$A_0(x) = A_1(x) = A_2(x) = A_3(x) = 1,$

$A_4(x) = 1 + x,$

$A_5(x) = 1 + 3x + x^2 = (1+x)^2 + x,$

$A_6(x) = 1 + 10x + 20x^2 + 10x^3 + x^4 = (1+x)^4 + 6x(1+x)^2 + 2x^2,$

$A_7(x) = 1 + 22x + 113x^2 + 190x^3 + 113x^4 + 22x^5 + x^6 = (1+x)^6 + 16x(1+x)^4 + 34x^2(1+x)^2 + 6x^3.$

Here again The Online Encyclopedia of Integer Sequences was very helpful.
OEIS [10], A245173, and its references to [1] and [9] showed that

$$A_{2k}(x) = \prod_{j=0}^{k-1} \frac{j!}{(k+j)!} \frac{(1-x)^{k^2+1}}{x^{k-1}} F_k^k(x) \tag{30}$$

and

$$A_{2k+1}(x) = \prod_{j=0}^{k} \frac{j!}{(k+j)!} \frac{(1-x)^{k^2+k+1}}{x^{k-1}} F_k^{k+1}(x). \tag{31}$$

The polynomials $F_a^b(x)$ are defined as

$$F_a^b(x) := \left(X^{a-1}D^a\right)^b \frac{1}{1-x}, \tag{32}$$

where $X$ denotes the multiplication operator $Xf(x) = xf(x)$ and $D$ denotes the differentiation operator $Df(x) = \frac{df(x)}{dx}.$

To give a sketch of the proof let

$$v_{a,b}(x) = \prod_{j=1}^{a}\prod_{i=1}^{b}(x+j-i) = \prod_{j=1}^{a}\prod_{i=1}^{b}(x+j+i-b-1). \tag{33}$$

**Proposition 3**

$$V_{a,b}(x) := \sum_{n\ge 0} v_{a,b}(n) x^{n+a-b-1} = F_a^b(x). \tag{34}$$

**Proof**

For $b \ge 1$ we get

$$V_{a,b}(x) = \sum_{n\ge 0} v_{a,b}(n)x^{n+a-b-1} = \sum_{n\ge 0}\prod_{j=1}^{a}\prod_{i=1}^{b}(n-i+j)x^{n+a-b-1} = \sum_{n\ge 0}\prod_{j=1}^{a}\prod_{i=1}^{b-1}(n-i+j)\prod_{j=1}^{a}(n-b+j)x^{n+a-b-1}$$

$$= x^{a-1}\sum_{n\ge 0}\prod_{j=1}^{a}\prod_{i=1}^{b-1}(n-i+j)D^a x^{n-b+a} = x^{a-1}D^a\sum_{n\ge 0}\prod_{j=1}^{a}\prod_{i=1}^{b-1}(n-i+j)x^{n+a-(b-1)-1} = x^{a-1}D^aV_{a,b-1}(x).$$

Therefore $V_{a,b}(x)$ satisfies $\left(x^{a-1}D^a\right)^b V_{a,0}(x) = \left(x^{a-1}D^a\right)^b \frac{x^{a-1}}{1-x}.$

Since $D^a\left(\frac{x^{a-1}}{1-x} - \frac{1}{1-x}\right) = D^a\left(1 + x + \cdots + x^{a-2}\right) = 0$ we get $V_{a,1}(x) = x^{a-1}D^a\frac{1}{1-x}.$

Therefore $\left(x^{a-1}D^a\right)^b \frac{x^{a-1}}{1-x} = \left(x^{a-1}D^a\right)^b \frac{1}{1-x}$ for $b \ge 1.$

By [1] and [9] we get

**Theorem 4**

$$V_{a,b}(x) := \sum_{n\geq 0} v_{a,b}(n)x^{n+a-b-1} = F_a^b(x) = v_{a,b}(b)\frac{x^{a-1}G_{a,b}(x)}{(1-x)^{ab+1}} = \prod_{j=0}^{b-1}\frac{(a+j)!}{j!}x^{a-1}\frac{G_{a,b}(x)}{(1-x)^{ab+1}} \tag{35}$$

*for* $b \geq 1$, *where* $F_a^b(x)$ *satisfies* $F_a^b(x) = x^{a-1}D^a F_a^{b-1}(x)$ *with* $F_a^0(x) = \dfrac{1}{1-x}$ *and where* $G_{a,b}(x)$ *is a gamma-positive polynomial of degree* $(a-1)(b-1)$.

We can write $r_k(n)$ as

$$r_k(n) = \frac{v_{\left\lfloor \frac{k}{2}\right\rfloor,\left\lfloor \frac{k+1}{2}\right\rfloor}\left(n+\left\lfloor \frac{k+1}{2}\right\rfloor\right)}{v_{\left\lfloor \frac{k}{2}\right\rfloor,\left\lfloor \frac{k+1}{2}\right\rfloor}\left(\left\lfloor \frac{k+1}{2}\right\rfloor\right)}, \tag{36}$$

because

$$\prod_{j=1}^{\left\lfloor \frac{k}{2}\right\rfloor}\prod_{i=1}^{\left\lfloor \frac{k+1}{2}\right\rfloor}\left(i+j+n-1\right)) = \prod_{j=1}^{\left\lfloor \frac{k}{2}\right\rfloor}\prod_{i=1}^{\left\lfloor \frac{k+1}{2}\right\rfloor}\left(i+j+n+\left\lfloor \frac{k+1}{2}\right\rfloor-\left\lfloor \frac{k+1}{2}\right\rfloor-1\right) = v_{\left\lfloor \frac{k}{2}\right\rfloor,\left\lfloor \frac{k+1}{2}\right\rfloor}\left(n+\left\lfloor \frac{k+1}{2}\right\rfloor\right)$$

This implies

**Corollary 5**

$$\sum_{n\geq 0} D_{2k}(n)x^n = \sum_{n\geq 0} r_{2k}(n)x^n = \frac{G_{k,k}(x)}{(1-x)^{k^2+1}}. \tag{37}$$

$$\sum_{n\geq 0} r_{2k+1}(n)x^n = \frac{G_{k,k+1}(x)}{(1-x)^{k^2+k+1}}. \tag{38}$$

**Proof**

We have

$v_{k,k}(n) = \prod_{j=1}^{k}\prod_{i=1}^{k}\left(n+j-i\right) = 0$ for $n < k$ because $n+j-k=0$ for some $j$.

$v_{k,k}(k) = \prod_{j=1}^{k}\prod_{i=1}^{k}\left(k+j-i\right) = \prod_{j=0}^{k-1}\prod_{i=0}^{k-1}\left(k+j-i\right) = \prod_{j=0}^{k-1}\frac{(k+j)!}{j!}.$

Therefore,

$$\sum_{n\geq 0}\sum_{n\geq 0} r_{2k}(n)x^n = \sum_{n\geq 0}\frac{v_{k,k}(n+k)}{v_{k,k}(k)}x^n = \prod_{j=0}^{k-1}\frac{j!}{(k+j)!}\sum_{n\geq 0} v_{k,k}(n+k)x^n$$

$$= \prod_{j=0}^{k-1}\frac{j!}{(k+j)!}\sum_{n\geq 0} v_{k,k}(n)x^{n-k}.$$

By Theorem 4 we have $\sum_{n\geq 0} v_{k,k}(n)x^{n-k} = x^{1-k}\sum_{n\geq 0} v_{k,k}(n)x^{n+k-k-1} = x^{1-k}\prod_{j=0}^{k-1}\frac{(k+j)}{j!}x^{k-1}\frac{G_{k,k}(x)}{(1-x)^{k^2+1}}.$

Therefore, we get (37) and in a similar way (38).

**Theorem 6**

*The generating function of* $D_{2k+1}(n)$ *satisfies*

$$\sum_{n\geq 0} D_{2k+1}(n)x^n = \frac{B_{2k+1}(x)}{\left(1+x^2\right)^{k^2+k+1}}. \tag{39}$$

*where* $B_{2k+1}(x)$ *is a polynomial of degree* $2k^2+1$ *which satisfies*

$$x^{2k^2+1} B_{2k+1}\left(\frac{1}{x}\right) = (-1)^k B_{2k+1}(x). \tag{40}$$

**Proof**

We have $B_{2k+1}(x) = \left(1+x^2\right)^{k^2+k+1} \sum_{n\geq 0} (-1)^{\binom{n}{2}} r_{2k+1}(n)x^n = B^0_{2k+1}(ix) + xB^1_{2k+1}(ix)$ with

$B^0_{2k+1}(x) = \left(1-x^2\right)^{k^2+k+1} \sum_{n\geq 0} r_{2k+1}(2n)x^{2n}$ and $B^1_{2k+1}(x) = \left(1-x^2\right)^{k^2+k+1} \sum_{n\geq 0} r_{2k+1}(2n+1)x^{2n}$.

By (29) we have $\sum_{n\geq 0} r_{2k+1}(n)x^n = \dfrac{A_{2k+1}(x)}{(1-x)^{k^2+k+1}}$. Therefore we get

$$B^0_{2k+1}(x) = \frac{1}{2}\left( A_{2k+1}(x)(1+x)^{k^2+k+1} + A_{2k+1}(-x)(1-x)^{k^2+k+1} \right)$$

and

$$B^1_{2k+1}(x) = \frac{1}{2x}\left( A_{2k+1}(x)(1+x)^{k^2+k+1} - A_{2k+1}(-x)(1-x)^{k^2+k+1} \right).$$

Then $\deg B^0_{2k+1}(x) = \deg B^1_{2k+1}(x) = 2k^2$.

For both have degree $\leq (k^2-k)+k^2+k+1 = 2k^2+1$. For $B^0_{2k+1}(x)$ the highest coefficient vanishes and therefore the degree is $2k^2$.

Thus we get

$$\sum_{n\geq 0} D_{2k+1}(2n)x^{2n} = \frac{B^0_{2k+1}(ix)}{\left(1+x^2\right)^{k^2+k+1}} \tag{41}$$

and

$$\sum_{n\geq 0} D_{2k+1}(2n+1)x^{2n+1} = \frac{xB^1_{2k+1}(ix)}{\left(1+x^2\right)^{k^2+k+1}}. \tag{42}$$

and

$$\sum_{n\geq 0} D_{2k+1}(n)x^n = \frac{B_{2k+1}(x)}{\left(1+x^2\right)^{k^2+k+1}} = \frac{B^0_{2k+1}(ix) + xB^1_{2k+1}(ix)}{\left(1+x^2\right)^{k^2+k+1}}. \tag{43}$$

with $\deg B^0_{2k+1}(ix) = \deg B^1_{2k+1}(ix) = 2k^2$.

For example we get

$$\sum_{n\geq 0} D_1(2n)x^{2n} = \sum_{n\geq 0} D_1(2n+1)x^{2n} = \sum_{n\geq 0} (-1)^n x^{2n} = \frac{1}{1+x^2} \quad \text{and} \quad \sum_{n\geq 0} D_1(n)x^n = \frac{1+x}{1+x^2}.$$

From $A_3(x)=1$ we get $B_3^0(x)=\frac{1}{2}\left((1+x)^3+(1-x)^3\right)=1+3x^2$ and

$B_3^1(x)=\frac{1}{2x}\left((1+x)^3-(1-x)^3\right)=3+x^2.$

This implies $B_3(x)=B_3^0(ix)+xB_3^1(ix)=1-3x^2+x\left(3-x^2\right)=1+3x-3x^2-x^3.$

From $A_5(x)=1+3x+x^2$ we get $2B_5^0(x)=\left(1+3x+x^2\right)(1+x)^7+\left(1-3x+x^2\right)(1-x)^7$

and $2xB_5^1(x)=\left(1+3x+x^2\right)(1+x)^7-\left(1-3x+x^2\right)(1-x)^7$ which gives

$B_5^0(ix)=1-43x^2+161x^4-105x^6+10x^8$ and $B_5^1(ix)=10-105x^2+161x^4-43x^6+x^8$ and

$$\sum_{n\geq 0}D_5(n)x^n=\frac{1+10x-43x^2-105x^3+161x^4+161x^5-105x^6-43x^7+10x^8+x^9}{\left(1+x^2\right)^7}.$$

To prove (40) we need

**Proposition 7**

*The polynomials* $B_{2k+1}^j(x)$ *satisfy*

$$B_{2k+1}^1(x)=x^{2k^2}B_{2k+1}^0\left(\frac{1}{x}\right). \tag{44}$$

**Proof**

$$x^{2k^2}B_{2k+1}^0\left(\frac{1}{x}\right)=\frac{1}{2}\left(x^{2k^2}A_{2k+1}\left(\frac{1}{x}\right)\left(\frac{x+1}{x}\right)^{k^2+k+1}+x^{2k^2}A_{2k+1}\left(-\frac{1}{x}\right)\left(\frac{x-1}{x}\right)^{k^2+k+1}\right)$$

$$=\frac{1}{2}\left(A_{2k+1}(x)\frac{1}{x}(x+1)^{k^2+k+1}-A_{2k+1}(-x)\frac{1}{x}(1-x)^{k^2k+1}\right)=B_{2k+1}^1(x).$$

Now we can prove (40): By (44) we get

$$B_{2k+1}(x)=B_{2k+1}^0(ix)+xB_{2k+1}^1(ix)=B_{2k+1}^0(ix)+x(ix)^{2k^2}B_{2k+1}^0\left(\frac{1}{ix}\right)=B_{2k+1}^0(ix)+(-1)^k x^{2k^2+1}B_{2k+1}^0\left(\frac{1}{ix}\right).$$

Therefore

$$x^{2k^2+1}B_{2k+1}\left(\frac{1}{x}\right)=x^{2k^2+1}\left(B_{2k+1}^0\left(\frac{1}{-ix}\right)+(-1)^k x^{-(2k^2+1)}B_{2k+1}^0\left(-ix\right)\right)$$

$$=(-1)^k B_{2k+1}^0\left(-ix\right)+x^{2k^2+1}B_{2k+1}^0\left(\frac{1}{-ix}\right)=(-1)^k B_{2k+1}^0\left(-ix\right)+(-1)^k x(-ix)^{2k^2}B_{2k+1}^0\left(\frac{1}{-ix}\right)$$

$$=(-1)^k\left(B_{2k+1}^0\left(-ix\right)+xB_{2k+1}^1(-ix)\right)=(-1)^k B_{2k+1}(x),$$

because $B_{2k+1}^j(x)=\sum_{\ell}c_\ell x^{2\ell}$.

The first terms are $B_1(x)=1,\ B_3(x)=1+3x-3x^2-x^3,$

$B_5(x)=1+10x-43x^2-105x^3+161x^4+161x^5-105x^6-43x^7+10x^8+x^9.$

As remarked in [9] $G_{a,b}(x)$ has a combinatorial interpretation as

$$G_{a,b}(x)=\sum_{T\in SYT\left(a^b\right)} x^{des(T)-b+1}. \tag{45}$$

For the polynomials $a_n(t),\ b_n(t)$ and $c_n(t)$ we have only partial results and some conjectures.

## 4. Results and conjectures for the Hankel determinants of $b_n(t)$.

By 2.2 we know that $b_n(t)$ are the moments of the orthogonal polynomials corresponding to $s_0=1,\ s_k=0$ for $k>0,\ t_{2j}=t$ and $t_{2j+1}=1.$

Let

$$D_k(n,t)=\det\left(b_{k+i+j}(t)\right)_{i,j=0}^{n-1} \tag{46}$$

with $D_k(0,t)=1.$

For $k=0$ we get from (18)

$$D_0(n,t)=t^{\left\lfloor \frac{n^2}{4}\right\rfloor}. \tag{47}$$

The orthogonal polynomials satisfy $p_n(x)=xp_{n-1}(x)-t_{n-2}p_{n-2}(x)$ for $n>1$ and $p_1(x)=x-1.$ Therefore, $p_{2n}(0)=(-t)^n,\ p_{2n-1}(0)=(-1)^n.$

If we write $d_k(n,t)=\dfrac{D_k(n,t)}{D_0(n,t)}$ and set $d(k,0,t)=1$ then we get from (18)

$$\begin{aligned} d_1(2n,t)&=(-t)^n,\\ d_1(2n+1,t)&=(-1)^n. \end{aligned} \tag{48}$$

I could not find a closed formula for $d_k(n,t).$ But using the condensation method we can derive closed expressions for small $k.$ (22) implies

$$t^{\left\lfloor \frac{n^2}{4}\right\rfloor+\left\lfloor \frac{(n-2)^2}{4}\right\rfloor}d_{k+2}(n-2,t)d_k(n,t)-t^{2\left\lfloor \frac{(n-1)^2}{4}\right\rfloor}d_{k+2}(n-1,t)d_k(n-1,t)+t^{2\left\lfloor \frac{(n-1)^2}{4}\right\rfloor}d_{k+1}(n-1,t)^2=0.$$

For $n\to 2n$ this gives

$$td_{k+2}(2n-2,t)d_k(2n,t)-d_{k+2}(2n-1,t)d_k(2n-1,t)+d_{k+1}(2n-1,t)^2=0 \tag{49}$$

and for $n\to 2n+1$

$$d_{k+2}(2n-1,t)d_k(2n+1,t)-d_{k+2}(2n,t)d_k(2n,t)+d_{k+1}(2n,t)^2=0. \tag{50}$$

Using these formulas, we get

$$d_2(n,t)=\sum_{j=0}^{n}t^j=\frac{1-t^{n+1}}{1-t}, \tag{51}$$

$$d_3(2n,t)=(-t)^n\sum_{j=0}^{2n}(2n+1-j)t^j=(-t)^n\frac{(2n+1)-(2n+2)t+t^{2n+2}}{(1-t)^2},$$
$$d_3(2n+1,t)=(-1)^n\sum_{j=0}^{2n+1}(j+1)t^j=(-1)^n\frac{1-\left(2n+3-(2n+2)t\right)t^{2n+3}}{(1-t)^2}, \tag{52}$$

and

$$d_4(n,t)=\binom{n+3}{3}t^n+\sum_{j=0}^{n-1}\binom{j+3}{3}\left(t^j+t^{2n-j}\right)=\frac{1-\left((n+2)^2-2\left(n^2+4n+3\right)t+(n+2)^2t^2\right)t^{n+1}+t^{2n+4}}{(1-t)^4}. \tag{53}$$

The first few polynomials are

$$\left(d_2(n,t)\right)_{n\geq 0}=\left(1,1+t,1+t+t^2,1+t+t^2+t^3,\cdots\right),$$
$$\left(d_4(n,t)\right)_{n\geq 0}=\left(1,1+4t+t^2,1+4t+10t^2+4t^3+t^4,1+4t+10t^2+20t^3+10t^4+4t^5+t^6,\cdots\right),$$
$$\left(d_6(n,t)\right)_{n\geq 0}=\begin{pmatrix}1,1+9t+9t^2+t^3,1+9t+45t^2+65t^3+45t^4+9t^5+t^6,\\ 1+9t+45t^2+165t^3+270t^4+270t^5+165t^6+45t^7+9t^8+t^9,\cdots\end{pmatrix}.$$

$$\left(d_3(n,t)\right)_{n\geq 0}=\left(1,1+2t,-t\left(3+2t+t^2\right),-\left(1+2t+3t^2+4t^3\right),t^2\left(5+4t+3t^2+2t^3+t^4\right)\cdots\right),$$

$$\left(d_5(n,t)\right)_{n\geq 0}=\left(1,1+6t+3t^2,-t\left(6+16t+21t^2+6t^3+t^4\right),-\left(1+6t+21t^2+56t^3+51t^4+30t^5+10t^6\right),\cdots\right).$$

Inspections of these polynomials suggest that

$$d_{2k}(n,t)=\sum_{j=0}^{n}\binom{j+k^2-1}{j}t^j+\cdots,$$
$$d_{2k+1}(2n+1,t)=(-1)^n\sum_{j=0}^{2n+1}\binom{j+k^2+k-1}{j}t^j+\cdots. \tag{54}$$

More precisely, define $r_k(n,t)$ by

$$r_{2k}(n,t)=d_{2k}(n,t),$$
$$r_{2k+1}(2n,t)=(-1)^n\frac{d_{2k+1}(2n,t)}{t^n}, \tag{55}$$
$$r_{2k+1}(2n+1,t)=(-1)^n d_{2k+1}(2n+1,t),$$

then we get

**Conjecture 8**

*The polynomials $r_{2k}(n,t)$ are palindromic and unimodal polynomials in $t$ with positive coefficients and degree $\deg_t r_{2k}(n,t)=kn$. They can be written in the form*

$$r_{2k}(n,t)=\frac{\sum_{j=0}^{k}(-1)^j c_{2k,j}(n,t)t^{j(n+j)}}{(1-t)^{k^2}} \tag{56}$$

*where* $c_{2k,j}(n,t)$ *are polynomials with* $\deg_t c_{2k,j}(n,t)=2j(k-j),$ *which satisfy*

$$c_{2k,k-j}(n,t)=t^{2j(k-j)}c_{2k,j}\left(n,\frac{1}{t}\right). \tag{57}$$

**Conjecture 9**

*The polynomials* $r_{2k+1}(n,t)$ *have positive coefficients and degree* $\deg_t r_{2k+1}(n,t)=kn.$ *They satisfy*

$$r_{2k+1}(2n,t)=\frac{\sum_{j=0}^{k}c_{0,2k+1,j}(n,t)t^{j(2n+j+1)}}{(1-t)^{k^2+k}} \tag{58}$$

*and*

$$r_{2k+1}(2n+1,t)=\frac{\sum_{j=0}^{k}c_{1,2k+1,j}(n,t)t^{j(2n+1+j)}}{(1-t)^{k^2+k}} \tag{59}$$

*where* $c_{0,2k+1,j}(n,t)$ *are polynomials with* $\deg_t c_{0,2k+1,j}(n,t)=(2j+1)(k-j)$ *and* $c_{1,2k+1,j}(n,t)$ *are polynomials with* $\deg_t c_{1,2k+1,j}(n,t)=(2k-2j+1)j.$

Consider for example $r_5(2n,t).$ Here we get

$$r(5,2n,t)=\frac{c_{1,5,0}(n,t)+t^{2n+2}c_{1,5,1}(n,t)+t^{2n+6}c_{1,5,2}(n,t)}{(1-t)^6} \tag{60}$$

with

$c_{1,5,0}(n,t)=(n+1)(2n+1)-2(n+1)(2n+3)t+(n+2)(2n+3)t^2,$

$c_{1,5,1}(n,t)=(n+1)(2n+3)^2-2(n+1)(n+2)(6n+5)t+(2n+1)(2n+3)(3n+5)t^2-2(n+1)^2(2n+3)t^3,$

$c_{1,5,2}(n,t)=1.$

For $n=0$ this reduces to

$$1=r_5(0,t)=\frac{(1-6t+6t^2)+\left(9-20t+15t^2-6t^3\right)t^2+t^6}{(1-t)^6}$$

and for $n=1$ to

$$6+16t+21t^2+6t^3+t^4=r_5(2,t)=\frac{\left(6-20t+15t^2\right)+\left(50-132t+120t^2-40t^3\right)t^{2+2}+t^{4+6}}{(1-t)^6}.$$

For the generating functions we get

**Conjecture 10**

$$\sum_{n\ge 0}d_{2k}(n,t)x^n=\frac{A_{2k}(x,t)}{\prod_{j=0}^{k}\left(1-t^jx\right)^{1+2j(k-j)}} \tag{61}$$

*where* $A_{2k}(x,t)$ *is a polynomial in* $x$ *with degree*

$$\deg A_{2k}(x,t)=\frac{(k-1)\left(k^2+k-3\right)}{3}=(k-1)^2+2\binom{k}{3}=\deg A_{2k}(x)+2\binom{k}{3}\ \textit{for}\ k\ge 1.$$

*It satisfies*

$$A_{2k}\left(\frac{1}{t^k x},t\right)\left(t^k x^2\right)^{\frac{(k-1)(k^2+k-3)}{6}}=A_{2k}(x,t). \tag{62}$$

Note that $\sum_{j=0}^{k}\left(1+2j(k-j)\right)=\frac{1}{3}(1+k)(k^2-k+3)=\left(k^2+1\right)+2\binom{k}{3}.$

Identity (62) is some analog of the palindromicity of $A_{2k}(x)$.

For example we have $A_0(x,t)=A_2(x,t)=1,\ \ A_4(x,t)=1+tx,$

$A_6(x,t)=1+4t(1+t)x+t^2\left(1-t+t^2\right)x^2-10t^4(1+t)x^3+t^5\left(1-t+t^2\right)x^4+4t^7(1+t)x^5+t^9x^6.$

Note that $A_6(x,1)=(1-x)^2A_6(x)$.

Some special cases are

$$\sum_{n\ge 0}d(2,n,t)x^n=\frac{1}{\left(1-x\right)\left(1-tx\right)}=\frac{1}{1-t}\frac{1}{1-x}-\frac{t}{1-t}\frac{1}{1-tx}$$

which give $d_2(n,t)=\frac{1-t^{n+1}}{1-t}=1+t+\cdots+t^n$,

$$\sum_{n\ge 0}d(4,n,t)x^n=\frac{1+tx}{\left(1-x\right)\left(1-tx\right)^3\left(1-t^2x\right)}$$

$$=\frac{1}{(1-t)^4}\frac{1}{1-x}-\frac{t\left(1+t^2\right)}{(1-t)^4}\frac{1}{1-tx}+\frac{t}{(1-t)^2}\frac{1}{(1-tx)^2}-\frac{2t}{(1-t)^2}\frac{1}{(1-tx)^3}+\frac{t^4}{(1-t)^4}\frac{1}{(1-t^2x)}$$

which gives

$$d(4,n,t)=\frac{1-(n+2)^2t^{n+1}+2\left(n^2+4n+3\right)t^{n+2}-(n+2)^2t^{n+3}+t^{2n+4}}{(1-t)^4}.$$

**Conjecture 11**

$$\sum_{n\ge 0}d_{2k+1}(2n,t)x^n=\frac{B_{0,k}(x,t)}{\prod_{j=0}^{k}\left(1+t^{2j+1}x\right)^{1+(2j+1)(k-j)}} \tag{63}$$

$$\sum_{n\ge 0}d_{2k+1}(2n+1,t)x^n=\frac{B_{1,k}(x,t)}{\prod_{j=0}^{k}\left(1+t^{2j}x\right)^{1+j(2(k-j)+1)}} \tag{64}$$

*with* $\deg B_{0,k}(x,t)=\deg B_{1,k}(x,t)=\frac{k(k+1)(2k+1)}{6}$ *and*

$$B_{1,k}(x,t)=(-1)^k x^{\frac{k(k+1)(2k+1)}{6}} t^{\frac{k(k+1)(k^2+k+1)}{3}} B_{0,k}\left(\frac{1}{t^{2k+1}x},t\right). \tag{65}$$

Note that

$$\frac{k(k+1)(2k+1)}{6}=k^2+\frac{k(k-1)(2k-1)}{6}=\deg B_{0,k}(x)+\frac{k(k-1)(2k-1)}{6}$$

and that the degree of the denominator is

$$k+1+\binom{k+1}{2}+2\binom{k+1}{3}=\left(k^2+k+1\right)+\frac{k(k-1)(2k-1)}{6}.$$

By (63) and(64) we get

**Corollary 12**

$$\sum_{n\geq 0} d_{2k+1}(n,t)x^n=\frac{B_k(x,t)}{\prod_{j=0}^{2k+1}\left(1+t^j x^2\right)^{jk-\frac{(j+1)(j-2)}{2}}} \tag{66}$$

with $\deg B_k(x,t)=\dfrac{k(k+1)(2k+1)}{3}+2(k+1)+2\dbinom{k+1}{2}+4\dbinom{k+1}{3}+1.$

Some special cases are

$$\sum_{n\geq 0} d(1,2n,t)x^n=\frac{1}{1+tx},$$

$$\sum_{n\geq 0} d(3,2n,t)x^n=\frac{1-t(1+2t)x}{(1+tx)^2\left(1+t^3x\right)},$$

$$\sum_{n\geq 0} d(5,2n,t)x^n=\frac{B_{0,2}(x,t)}{(1+tx)^3\left(1+t^3x\right)^4\left(1+t^5x\right)}$$

with

$$B_{0,2}(x,t)=1-t\left(3+16t+17t^2+6t^3\right)x+2t^4\left(8+22t+21t^2+8t^3\right)x^2+t^5\left(9+34t+28t^2-13t^4-2t^5\right)x^3$$
$$-t^8\left(14+32t+36t^2+12t^3+t^4\right)x^4+t^{11}\left(1+6t+3t^2\right)x^5.$$

For $t=1$ this reduces to

$$1-42x+118x^2+56x^3-95x^4+10x^5=(1+x)\left(1-43x+161x^2-105x^3+10x^4\right)=(1+x)B_{0,2}(x).$$

## 5. Hankel determinants of $a_n(t)$ and their generating functions.

Let

$$\Delta_k(n,t)=\det\left(a_{k+i+j}(t)\right)_{i,j=0}^{n-1} \tag{67}$$

with $\Delta_k(0,t)=1.$

The polynomials $a_n(t)$ are the moments of the orthogonal polynomials with $s_0=1,$ $s_k=0$ for $k>0$ and $t_k=t.$

Therefore by (18) we get

$$\Delta_0(n,t) = t^{\binom{n}{2}}. \tag{68}$$

Let

$$\delta_k(n,t) = \frac{\Delta_k(n,t)}{\Delta_0(n,t)} = \frac{\Delta_k(n,t)}{t^{\binom{n}{2}}}. \tag{69}$$

Again by (18) we get

$$\delta_1(n,t) = (-1)^{\binom{n}{2}} t^{\left\lfloor \frac{n}{2} \right\rfloor}. \tag{70}$$

**Conjecture 13**

*The polynomials* $\delta_k(n,t)$ *have degree* $\deg_t \delta_k(n,t) = \left\lfloor \frac{kn}{2} \right\rfloor$.

*More precisely*

$$\delta_k(n,t) = (-1)^{k\binom{n}{2}} t^{\left\lfloor \frac{k(n-1)+1}{2} \right\rfloor} \rho_k(n,t) \tag{71}$$

*where* $\rho_k(n,t)$ *is a polynomial with positive coefficients with* $\deg_t \rho_k(n,t) = \left\lfloor \frac{k}{2} \right\rfloor$.

Let us consider $\rho_k(n,t)$ for small $k$ and $n$.

$\rho_1(n,t) = 1,$

$\rho_2(n,t) = \left\lfloor \frac{n+1}{2} \right\rfloor + t\left\lfloor \frac{n+2}{2} \right\rfloor,$

$\rho_3(n,t) = \left\lfloor \frac{(n+1)^2}{4} \right\rfloor + t\left\lfloor \frac{(n+2)^2}{4} \right\rfloor,$

$\rho_4(2n+1,t) = \frac{(n+1)(n+2)(2n+3)}{6}(1+t)\left(n+1+(n+2)t\right),$

$\rho_4(2n,t) = \frac{(n+1)^2}{6}\left(n(2n+1)+4n(n+2)t+(n+2)(2n+3)t^2\right).$

For the generating functions we get

**Conjecture 14**

$$\sum_{n\geq 0} \delta_{2k}(n,t)x^n = \frac{\mathrm{A}_{2k}(x,t)}{\left(1-t^k x\right)\left(1-t^{2k}x^2\right)^{k^2}} \tag{72}$$

*with* $\deg \mathrm{A}_{2k}(x,t) = 2k(k-1)+1$

*and*

$$\sum_{n\ge 0}\delta_{2k+1}(n,t)x^n=\frac{\mathrm{A}_{2k+1}(x,t)}{\left(1+t^{2k+1}x^2\right)^{k^2+k+1}} \tag{73}$$

*with* $\deg \mathrm{A}_{2k+1}(x,t)=2k^2+1.$

The first few $\mathrm{A}_k(x,t)$ are

$\mathrm{A}_1(x,t)=1+x,\ \mathrm{A}_2(x,t)=1+x,\ \mathrm{A}_3(x,t)=1+(1+2t)x-t^2(2+t)x^2-t^3x^3,$

$\mathrm{A}_4(x,t)=1+\left(1+3t+t^2\right)x+t^2\left(1+t\right)\left(1+4t\right)x^2+t^4\left(1+t\right)\left(4+t\right)x^3+t^6\left(1+3t+t^2\right)x^4+t^8x^5.$

It seems that the following analogs of palindromicity hold:

**Conjecture 15**

$$x^{2k(k-1)+1}t^{2k\left\lfloor\frac{(2k-1)^2}{4}\right\rfloor}\mathrm{A}_{2k}\left(\frac{1}{x},\frac{1}{t}\right)=x^{2k(k-1)+1}t^{2k^2(k-1)}\mathrm{A}_{2k}\left(\frac{1}{x},\frac{1}{t}\right)=\mathrm{A}_{2k}\left(x,t\right), \tag{74}$$

$$(-1)^k x^{2k^2+1}t^{(2k+1)k^2}\mathrm{A}_{2k+1}\left(\frac{1}{x},\frac{1}{t}\right)=\mathrm{A}_{2k+1}(x,t). \tag{75}$$

## 6. Hankel determinants of $c_n(t)$ and their generating functions

Let

$$\mathbf{D}_k(n,t)=\det\left(c_{k+i+j}(t)\right)_{i,j=0}^{n-1} \tag{76}$$

with $\mathbf{D}_k(0,t)=1.$

Then

$$\mathbf{D}_0(n,t)=t^{\binom{n}{2}}. \tag{77}$$

**Conjecture 16**

*The polynomial*

$$\mathbf{d}_k(n,t)=\frac{\mathbf{D}_k(n,t)}{t^{\binom{n}{2}}} \tag{78}$$

*has degree* $\deg_t \mathbf{d}_k(n,t) = (k-1)n$ *and* $(-1)^{k\binom{n}{2}}\mathbf{d}_k(n,t)$ *is palindromic with positive coefficients. Moreover*

$$[t^j]\mathbf{d}_{2k}(n,t) = [t^j]\frac{1}{(1-t)^k\left(1-t^2\right)^{k^2-k}} \tag{79}$$

*and*

$$[t^j](-1)^{\binom{n}{2}}\mathbf{d}_{2k+1}(n,t) = [t^j]\frac{1}{(1-t)^k\left(1-t^2\right)^{k^2}} \tag{80}$$

*for* $0 \le j \le n.$

We get $(-1)^{\binom{n}{2}}\mathbf{d}_1(n,t) = 1,\ \mathbf{d}_2(n,t) = 1 + t + \cdots + t^n,$
$(-1)^{\binom{n}{2}}\mathbf{d}_3(n,t) = \begin{bmatrix} n+2 \\ 2 \end{bmatrix}_t = \dfrac{\left(1-t^{n+1}\right)\left(1-t^{n+2}\right)}{\left(1-t\right)\left(1-t^2\right)},$

We have $\dfrac{1}{(1-t)^2\left(1-t^2\right)^2} = 1 + 2t + 5t^2 + 8t^3 + 14t^4 + \cdots$ and the first few polynomials $\mathbf{d}_4(n,t)$ are $\mathbf{1},\ \mathbf{1} + \mathbf{2}t + 2t^2 + t^3,\ \mathbf{1} + \mathbf{2}t + \mathbf{5}t^2 + 4t^3 + 5t^4 + 2t^5 + t^6,$
$\mathbf{1} + \mathbf{2}t + \mathbf{5}t^2 + \mathbf{8}t^3 + 9t^4 + 9t^5 + 8t^6 + 5t^7 + 2t^8 + t^9, \cdots.$

Similarly we get $\dfrac{1}{(1-t)^2\left(1-t^2\right)^4} = 1 + 2t + 7t^2 + 12t^3 + 27t^4 + \cdots$ and the first few polynomials $\mathbf{d}_5(n,t)$ are $\mathbf{1},\ \mathbf{1} + \mathbf{2}t + 4t^2 + 2t^3 + t^4,\ -\left(\mathbf{1} + \mathbf{2}t + \mathbf{7}t^2 + 8t^3 + 14t^4 + 8t^5 + 7t^6 + 2t^7 + t^8\right),$

$-\left(\mathbf{1} + \mathbf{2}t + 7t^2 + \mathbf{12}t^3 + 22t^4 + 26t^5 + 35t^6 + 26t^7 + 22t^8 + 12t^9 + 7t^{10} + 2t^{11} + t^{12}\right), \cdots.$

More precisely we get

**Conjecture 17**

*The polynomials* $(-1)^{k\binom{n}{2}}\mathbf{d}_k(n,t)$ *can be written in the form*

$$\mathbf{d}_{2k}(n,t) = \frac{\sum_{j=0}^{2k-1}(-1)^j c_{2k,j}(n,t)t^{jn+\binom{j+1}{2}}}{(1-t)^k\left(1-t^2\right)^{k^2-k}} \tag{81}$$

*and*

$$(-1)^{\binom{n}{2}}\mathbf{d}_{2k+1}(n,t)=\frac{\sum_{j=0}^{2k}(-1)^j c_{2k+1,j}(n,t)t^{jn+\binom{j+1}{2}}}{(1-t)^k\left(1-t^2\right)^{k^2}} \tag{82}$$

*with* $\deg_t c_{k,j}(n,t)=j\left(k-1-j\right)$ *and satisfy*

$$c_{k,k-1-j}(n,t)=(-1)^{\left\lfloor \frac{k}{2}\right\rfloor-1}t^{j(k-1-j)}c_{k,j}\left(n,\frac{1}{t}\right). \tag{83}$$

For example we get

$$\mathbf{d}_2(n,t)=\frac{1-t^{n+1}}{1-t},$$

with $c_{2,0}(n,t)=c_{2,1}(n,t)=1,$

$$\mathbf{d}_4(n,t)=\frac{1-t^{n+1}\left(n+2+t-(n+2)t^2\right)+t^{2n+3}\left(n+2-t-(n+2)t^2\right)+t^{3n+6}}{(1-t)^2\left(1-t^2\right)^2}$$

with $c_{4,0}(n,t)=1,\ c_{4,3}(n,t)=-1,$

$$c_{4,1}(n,t)=n+2+t-(n+2)t^2,\ \ c_{4,2}(n,t)=n+2-t-(n+2)t^2=-t^2c_{4,1}\left(n,\frac{1}{t}\right).$$

For $\mathbf{d}_6(n,t)$ we get

$$c_{6,0}(n,t)=c_{6,5}(n,t)=1,$$

$$c_{6,1}(n,t)=\binom{n+3}{2}+(n+3)t-(n+2)(n+4)t^2-(n+3)t^3+\binom{n+4}{2}t^4$$

$$c_{6,4}(n,t)=\binom{n+4}{2}-(n+3)t-(n+2)(n+4)t^2+(n+3)t^3+\binom{n+4}{2}t^4=t^4c_{6,1}(n,t),$$

$$c_{6,2}(n,t)=\binom{n+3}{2}(n+3)-\binom{n+4}{2}t-\binom{n+3}{2}(3n+11)t^2+\left((n+3)^2-2\right)t^3+\binom{n+4}{2}(3n+7)t^4$$
$$-\binom{n+3}{2}t^5-\binom{n+4}{2}(n+3)t^6$$

$$c_{6,3}(n,t)=-\binom{n+3}{2}(n+3)t^6-\binom{n+4}{2}t^5+\binom{n+3}{2}(3n+11)t^4+\left((n+3)^2-2\right)t^3-\binom{n+4}{2}(3n+7)t^2$$
$$-\binom{n+3}{2}t+\binom{n+4}{2}(n+3)=t^6c_{6,2}\left(n,\frac{1}{t}\right)$$

For $\mathbf{d}_3(n,t)$ we get $c_{3,0}(n,t)=c_{3,1}(n,t)=1,\ c_{3,1}(n,t)=1+t.$

For $\mathbf{d}_5(n,t)$ we get $c_{5,1}(n,t)=\left(1-t^2\right)\left(n+2+(n+3)t\right),$
$c_{5,2}(n,t)=(n+2)(n+3)-t-2(n+2)(n+3)t^2-t^3+(n+2)(n+3)t^4.$

For the generating functions we get

$$\sum_{n\geq 0}\mathbf{d}_1(n,t)x^n=\frac{1+x}{1+x^2},\ \sum_{n\geq 0}\mathbf{d}_2(n,t)x^n=\frac{1}{(1-x)(1-tx)},$$

$$\sum_{n\geq 0}\mathbf{d}_3(n,t)x^n=\frac{1+(1+t)^2x+t^2x^2}{\left(1+x^2\right)\left(1+t^2x^2\right)\left(1+t^4x^2\right)},$$

$$\sum_{n\geq 0}\mathbf{d}_4(n,t)x^n=\frac{1-t^3x^2}{(1-x)(1-tx)^2\left(1-t^2x\right)^2\left(1-t^3x\right)},$$

**Conjecture 18**

*The generating functions are*

$$\sum_{n\geq 0}\mathbf{d}_{2k}(n,t)x^n=\frac{\mathbf{A}_{2k}(x,t)}{\prod_{j=0}^{2k-1}\left(1-t^jx\right)^{kj-\frac{(j-1)(j+2)}{2}}} \tag{84}$$

*with*

$$\mathbf{A}_{2k}\left(\frac{1}{t^{2k-1}x},t\right)x^{\frac{(k-1)k(2k-1)}{3}}t^{\frac{(k-1)k(2k-1)^2}{6}}=\mathbf{A}_{2k}(x,t) \tag{85}$$

*and*

$$\sum_{n\geq 0}\mathbf{d}_{2k-1}(n,t)x^n=\frac{\mathbf{A}_{2k-1}(x,t)}{\prod_{j=0}^{2k-2}\left(1+t^{2j}x^2\right)^{kj-\frac{(j-1)(j+2)}{2}-\left\lfloor\frac{j+1}{2}\right\rfloor}} \tag{86}$$

*with*

$$(-1)^k\mathbf{A}_{2k+1}\left(\frac{1}{t^{2k}x},t\right)x^{\frac{2k(k-1)(2k-1)}{3}-2(k-1)^2+1}t^{k\left(\frac{2k(k-1)(2k-1)}{3}-2(k-1)^2+1\right)}=\mathbf{A}_{2k+1}(x,t). \tag{87}$$

Finally let us note that the Hankel matrices of the sequence $\left(c_n(-1)\right)_{n\geq 0}=\left(1,1,0,1,0,2,0,5,0,14,0,\cdots\right)$ are checkerboard matrices. Their determinants can be reduced to the Hankel determinants of the Catalan numbers.

$$\det\left(c(2k+i+j,-1)\right)_{i,j=0}^{2n-1}=(-1)^n\left(\det\left(C_{k+i+j}\right)_{i,j=0}^{n-1}\right)^2,$$

$$\det\left(c(2k+1+i+j,-1)\right)_{i,j=0}^{2n-1}=0 \text{ for } k>0,$$
$$\det\left(c(2k-1+i+j,-1)\right)_{i,j=0}^{2n-1}=\det\left(C_{k+i+j}\right)_{i,j=0}^{n-1}\det\left(C_{k+i+j+1}\right)_{i,j=0}^{n-1},$$
$$\det\left(c(2k-1+i+j,-1)\right)_{i,j=0}^{2n}=\det\left(C_{k+i+j}\right)_{i,j=0}^{n-1}\det\left(C_{k+i+j+1}\right)_{i,j=0}^{n}.$$

## 7. Conjectures for the Hankel determinants of the numbers $b^{(r)}(n)=\binom{n}{\left\lfloor \frac{n-r}{2}\right\rfloor}$.

In [4] it has been proved that many Hankel determinants of the sequences $\binom{2n+r}{n}$ show an interesting modular pattern. The same seems to be true for the Hankel determinants of the numbers $b^{(r)}(n)=\binom{n}{\left\lfloor \frac{n-r}{2}\right\rfloor}$.

Let

$$D_k^{(r)}(n)=\det\left(b^{(r)}(k+i+j)\right)_{i,j=0}^{n-1} \tag{88}$$

with $D_k^{(r)}(0)=1$.

For $r\in\mathbb{N}$ we have

$$\begin{aligned} &D_0^{(r)}((2r+1)n)=1,\\ &D_0^{(r)}((2r+1)n+r+1)=(-1)^{\binom{r+1}{2}} \end{aligned} \tag{89}$$

and $D_0^{(r)}((2r+1)n+i)=0$ else.

For $k=1$ we see that

$$D_1^{(r)}((4r+2)n+i)\neq 0 \quad \text{for } i\in\{0,2r+1\}\cup\{r,3r+1\}\cup\{r+1,2r\}.$$

and $D_1^{(r)}((4r+2)n+i)=0$ else.

More precisely we get

$$\begin{aligned} &D_1^{(r)}((4r+2)n)=D_1^{(r)}((4r+2)n+2r+1)=(-1)^n,\\ &D_1^{(r)}((4r+2)n+r)=D_1^{(r)}((4r+2)n+3r+1)=(-1)^{n+\binom{r}{2}},\\ &D_1^{(r)}((4r+2)n+r+1)=(-1)^{\binom{r}{2}}D_1^{(r)}\left((4r+2)n+2r\right)=2(-1)^n. \end{aligned} \tag{90}$$

For example

$$\left(D_1^{(1)}(n)\right)_{n\geq 0}=\left(1,1,2,1,1,0,-1,-1,-2,-1,-1,0,1,1,2,1,1,0,\cdots\right)$$

$$\left(D_1^{(2)}(n)\right)_{n\ge 0}=\left(1,0,-1,2,-2,1,0,-1,0,0,-1,0,1,-2,2,-1,0,1,0,0,\cdots\right).$$

For $k=2$ we get

$$D_2^{(r)}((4r+2)n+i)\neq 0 \quad \text{for } i\in\{0,2r+1\}\cup\{r-1,3r\}\cup\{r,2r\}\cup\{r+1,2r-1\}\cup\{3r+1,4r+1\}.$$

More precisely we get for $r\ge 2$

$$\begin{aligned} &D_2^{(r)}\left((4r+2)n\right)=D_2^{(r)}\left((4r+2)n+2r+1\right)=1,\\ &D_2^{(r)}\left((4r+2)n+r-1\right)=D_2^{(r)}\left((4r+2)n+3r\right)=(-1)^{\binom{r-1}{2}}, \end{aligned} \tag{91}$$

$$\begin{aligned} &D_2^{(r)}\left((4r+2)n+r\right)=(2r+1)(2n+1)-2\left\lfloor\frac{r+1}{2}\right\rfloor,\\ &D_2^{(r)}((4r+2)n+2r)=(-1)^{\binom{r+1}{2}}\left((2r+1)(2n+1)+2\left\lfloor\frac{r+1}{2}\right\rfloor\right), \end{aligned} \tag{92}$$

$$\begin{aligned} &D_2^{(r)}\left((4r+2)n+r+1\right)=4(-1)^{\binom{r+1}{2}},\\ &D_2^{(r)}\left((4r+2)n+2r-1\right)=-4, \end{aligned} \tag{93}$$

$$\begin{aligned} &D_2^{(r)}\left((4r+2)n+3r+1\right)=(4r+2)(n+1),\\ &D_2^{(r)}\left((4r+2)n+4r+1\right)=(-1)^{\binom{r+1}{2}}(4r+2)(n+1). \end{aligned} \tag{94}$$

For example

$$\left(D_2^{(2)}(n)\right)_{n\ge 0}=\left(1,1,3,-4,-7,1,1,10,0,-10,1,1,13,-4,-17,1,1,20,0,-20,\cdots\right),$$

$$\left(D_2^{(3)}(n)\right)_{n\ge 0}=\left(1,0,-1,3,4,-4,11,1,0,-1,14,0,0,14,1,0,-1,17,4,-4,25,1,0,-1,28,0,0,28,\cdots\right).$$

For $r=1$ we get

$$\left(D_2^{(1)}(n)\right)_{n\ge 0}=\left(1,1,-5,1,6,-6,1,7,11,1,12,-12,\cdots\right).$$

Here we have

$$\begin{aligned} &D_2^{(1)}(6n)=1,\ \ D_2^{(1)}(6n+3)=1,\\ &D_2^{(1)}(6n+1)=6n+1,\ \ D_2^{(1)}(6n+2)=-(6n+5),\\ &D_2^{(1)}(6n+4,1)=6(n+1),\ \ D_2^{(1)}(6n+5,1)=-6(n+1). \end{aligned} \tag{95}$$

For $r\le k$ the Hankel determinants $D_k^{(r)}((4r+2)n+i)$ seem to be of the form $(-1)^{kn}\,p_{k,i}(n)$ for some polynomial $p_{k,i}(x)$ with $\deg p_{2k,i}(x)\le k^2$ and $\deg p_{2k+1,i}(x)\le k^2+k.$ For example, for $D_4^{(r)}((4r+2)n+1)$ we get $D_4^{1}(6n+1)=(2n+1)(2n+3)(18n^2+21n+2)$,

$$D_4^{2}(10n+1)=\frac{(5n+2)\left(500n^3+700n^2+165n-6\right)}{3},\ D_4^{(3)}(14n+1)=\frac{(14n+3)\left(196n^2+63n-1\right)}{3}.$$

**The general case for $r > k$.**

For $r > k > 2$ the Hankel determinants $D_k^{(r)}\left((4r+2)n+i\right)$ vanish except for

$i \in \{0,-1,\cdots,-(k-1)\} \cup \{r+1,r,\cdots,r-(k-1)\} \cup \{2r+1,2r,\cdots,2r-(k-1)\} \cup \{3r+1,3r,\cdots,3r-(k-2)\}$.

For $i \in \{0,r-k+1\} \cup \{r+1,2r-k+1\} \cup \{2r+1,3r-k+2\} \cup \{3r+1,4r-k+3\}$ we get some nice evaluations.

For $r \geq k-1$ we get

$$D_k^{(r)}\left((4r+2)n\right) = D_k^{(r)}\left((4r+2)n+r-k+1\right)(-1)^{\binom{r-k+1}{2}} = (-1)^{kn}, \tag{96}$$

for $r \geq k$ we get

$$(-1)^{kn}(-1)^{(k+1)\binom{r+1}{2}} D_k^{(r)}\left((4r+2)n+r+1\right) = -D_k^{(r)}\left((4r+2)n+2r-k+1\right) = 2^k, \tag{97}$$

for $r \geq 1$ we get

$$D_k^{(r)}\left((4r+2)n+2r+1\right) = (-1)^{\binom{r-1}{2}} D_k^{(r)}\left((4r+2)n+3r-k+2,r\right) = 1, \tag{98}$$

for $r \geq k-2$ we get

$$\begin{aligned} D_k^{(r)}\left((4r+2)n+3r+1,r\right) &= (-1)^{k\left(n+\binom{r}{2}\right)}\left((4r+2)(n+1)\right)^{k-1}, \\ D_k^{(r)}\left((4r+2)n+4r-k+3,r\right) &= (-1)^{kn+(k+1)\binom{r+k-1}{2}+[k\equiv 1 \bmod 4]}\left((4r+2)(n+1)\right)^{k-1}. \end{aligned} \tag{99}$$

Further we get the following identities:

For $r \geq k-2$

$$D_k^{(r)}\left((4r+2)n+r\right)(-1)^{\binom{r+1-k}{2}} + D_k^{(r)}\left((4r+2)n+2(r+1)-k\right) = (-1)^{kn+(k+1)\binom{r+1}{2}} 2^k \binom{\left\lfloor \frac{r+k-1}{2} \right\rfloor}{k-1}, \tag{100}$$

for $r \geq k-3$ and $k \geq 2$

$$D_k^{(r)}\left((4r+2)n+2r\right)(-1)^{\binom{r+1-k}{2}} + D_k^{(r)}\left((4r+2)n+3r+3-k\right) = (-1)^{kn} c\left(k, \left\lfloor \frac{r+2-k}{2} \right\rfloor\right) \tag{101}$$

Here $c(k,n)$ is defined by $c(k,n) = \sum_{j=0}^{n} c(k-1,j)$ with $c(2,n) = 2n+1$.

For $r \ge k-3$ we get

$$D_k^{(r)}\left((4r+2)n+4(r+1)-k\right)-D_k^{(r)}\left((4r+2)n+3r\right)(-1)^{\binom{r+1-k}{2}}$$
$$=(-1)^{k\left(n+\binom{r}{2}\right)}\binom{\left\lfloor \frac{r+k-3}{2}\right\rfloor}{k-3}\left((4r+2)(n+1)\right)^{k-1}. \tag{102}$$

## 8. Some q-analogs

Let $[n]_q=1+q+\cdots+q^{n-1}=\dfrac{1-q^n}{1-q}$, $[n]_q!=[1]_q[2]_q\cdots[n]_q$, $\begin{bmatrix} n \\ k \end{bmatrix}_q=\dfrac{[n]_q!}{[k]_q![n-k]_q!}$ for $0\le k\le n$, $\begin{bmatrix} n \\ k \end{bmatrix}_q=0$ else and $(x;q)_n=\prod_{j=0}^{n-1}\left(1-q^jx\right)$.

I could not find polynomial $q-$analogs of $b(n)$ with nice Hankel determinants but the following rational functions give nice values.

Let

$$b(n,q)=\begin{bmatrix} n \\ \left\lfloor \frac{n}{2}\right\rfloor \end{bmatrix}_q \frac{2^n}{\left(-q;q\right)_{\left\lfloor \frac{n}{2}\right\rfloor}\left(-q;q\right)_{\left\lfloor \frac{n+1}{2}\right\rfloor}}. \tag{103}$$

These are moments of the orthogonal polynomials which correspond to the values

$$t_n=\frac{4q^{n+1}}{\left(1+q^{n+1}\right)^2}, s_0=\frac{2}{1+q}, s_n=\frac{2q^n(1-q)}{\left(1+q^n\right)\left(1+q^{n+1}\right)}. \tag{104}$$

This follows from the fact that (16) holds with

$$c(2n,2k,q)=\begin{bmatrix} 2n \\ n-k \end{bmatrix}_q 2^{2n-2k}\frac{\left(-q;q\right)_{2k}}{\left(-q;q\right)_{n-k}\left(-q;q\right)_{n+k}},$$

$$c(2n-1,2k,q)=\begin{bmatrix} 2n-1 \\ n-k-1 \end{bmatrix}_q 2^{2n-1-2k}\frac{\left(-q;q\right)_{2k}}{\left(-q;q\right)_{n-k-1}\left(-q;q\right)_{n+k}},$$

$$c(2n,2k-1,q)=\begin{bmatrix} 2n \\ n-k \end{bmatrix}_q 2^{2n-2k+1}\frac{\left(-q;q\right)_{2k-1}}{\left(-q;q\right)_{n-k}\left(-q;q\right)_{n+k}},$$

$$c(2n-1,2k-1,q)=\begin{bmatrix} 2n-1 \\ n-k \end{bmatrix}_q 2^{2n-2k}\frac{\left(-q;q\right)_{2k-1}}{\left(-q;q\right)_{n-k}\left(-q;q\right)_{n+k-1}}$$

and that $c(2n,0,q)=\begin{bmatrix}2n\\ n\end{bmatrix}_q 2^{2n}\dfrac{1}{(-q;q)_n(-q;q)_n}=b(2n,q)$ and

$$c(2n-1,0,q)=\begin{bmatrix}2n-1\\ n-1\end{bmatrix}_q 2^{2n-1}\frac{1}{(-q;q)_{n-1}(-q;q)_n}=b(2n-1,q).$$

Let $D_k(n,q)=\det\left(b(i+j+k,q)\right)_{i,j=0}^{n-1}$.

Then (18) gives

$$D_0(n,q)=\frac{2^{n(n-1)}q^{\binom{n+1}{3}}}{\prod_{j=0}^{n-1}(-q;q)_j^2}. \tag{105}$$

If we set

$$d_k(n,q)=\frac{D_k(n,q)}{D_0(n,q)} \tag{106}$$

and

$$r_k(n,q)=\prod_{i=1}^{\left\lfloor\frac{k}{2}\right\rfloor}\prod_{j=1}^{\left\lceil\frac{k}{2}\right\rceil}\frac{[i+j+n-1]_q}{[i+j-1]_q} \tag{107}$$

then we get

$$\begin{aligned}d_{2k}(2n,q)&=q^{2kn^2}\frac{2^{4kn}r_{2k}(2n)}{\prod_{j=1}^{k}(-q^j;q)_{2n}^2},\\ d_{2k}(2n+1,q)&=q^{2kn(n+1)}\frac{2^{2k(2n+1)}r_{2k}(2n+1)}{\prod_{j=1}^{k}(-q^j;q)_{2n+1}^2}\end{aligned} \tag{108}$$

and

$$\begin{aligned}d_{2k+1}(2n,q)&=(-1)^n q^{((2k+1)n+1)n}\frac{2^{2(2k+1)n}r_{2k+1}(2n)}{\prod_{j=1}^{k}(-q^j;q)_n^2\cdot(-q^{k+1};q)_{2n}},\\ d_{2k+1}(2n+1,q)&=(-1)^n q^{(2k+1)n(n+1)}\frac{2^{(2k+1)(2n+1)}r_{2k+1}(2n+1)}{\prod_{j=1}^{k}(-q^j;q)_{2n+1}^2\cdot(-q^{k+1};q)_{2n+1}}.\end{aligned} \tag{109}$$